\documentclass{amsart}
\usepackage{amsmath,amsthm,amssymb,amsfonts}
\usepackage{graphicx} 
\usepackage{epstopdf}

\makeatletter

\def\proof{\@ifstar{P\,r\,o\,o\,f}{P\,r\,o\,o\,f.\ }}

\renewcommand\th@remark{%
  \thm@headfont{\bfseries}%
  \normalfont 
  \thm@preskip\topsep \divide\thm@preskip\tw@
  \thm@postskip\thm@preskip
}

\renewenvironment{equation}{\refstepcounter{equation}$$}{\eqno{(\thesection.\theequation)}$$}

\makeatother

\newcounter{Example}[section]
\newcounter{Th}[section] \newcounter{Pr}[section] \newcounter{Lm}[section]
\newcounter{Remark}[section]
\newcounter{Def}[section]
\newcounter{lcounter}[section]
\newcounter{Corol}[section]
\newcounter{Conj}[section]

\newenvironment{Th}[1][\relax]
    {\medspace\refstepcounter{Th}T\,h\,e\,o\,r\,e\,m \arabic{section}.\theTh.\ \it}
    {\rm\medspace}

\newenvironment{Th.}[1][\relax]
    {\medspace\refstepcounter{Th}T\,h\,e\,o\,r\,e\,m \arabic{section}.\theTh.\ \it}
    {\rm\medspace}

\newenvironment{Pr.}[1][\relax]
    {\medspace\refstepcounter{Pr}P\,r\,o\,p\,o\,s\,i\,t\,i\,o\,n \arabic{section}.\thePr.\ \it}
    {\rm\medspace}

\newenvironment{Lm}[1][\relax]
    {\medspace\refstepcounter{Lm}L\,e\,m\,m\,a \arabic{section}.\theLm.\ \it}
    {\rm\medspace}

\newenvironment{Corol}[1][\relax]
    {\medspace\refstepcounter{Corol} C\,o\,r\,o\,l\,l\,a\,r\,y \arabic{section}.\theCorol.\ \it}
    {\rm\medspace}

\newenvironment{Conj}[1][\relax]
    {\medspace\refstepcounter{Conj} C\,o\,n\,j\,e\,c\,t\,u\,r\,e \arabic{section}.\theConj.\ \it}
    {\rm\medspace}

\newenvironment{Remark}[1][\relax]
    {\medspace\refstepcounter{Remark}R\,e\,m\,a\,r\,k \arabic{section}.\theRemark.\rm\ }
    {\medspace}

\newenvironment{Example}[1][\relax]
    {\medspace\refstepcounter{Example}E\,x\,a\,m\,p\,l\,e \arabic{section}.\theExample.\rm\ }
    {\medspace}

\newenvironment{Def}[1][\relax]
    {\medspace\refstepcounter{Def}D\,e\,f\,i\,n\,i\,t\,i\,o\,n \arabic{section}.\theDef.\rm\ }
    {\medspace}

\newenvironment{Def.}[1][\relax]
    {\medspace\refstepcounter{Def}D\,e\,f\,i\,n\,i\,t\,i\,o\,n \arabic{section}.\theDef.\rm\ }
    {\medspace}

\def\au#1{\emph{#1}}
\def\tit#1{{#1}}
\def\R{{\mathbb R}}  
\def\N{{\mathbb N}}  
\def\H{\mathcal{H}} 
\def\E{E} 
\def\cl{\mathop{\rm cl\,}}

\def\Bigcap{\bigcap\limits}

\def\co{\mbox{\rm co}\,}

\def\aff {\mbox{\rm aff}\,}

\def\ll{\lambda}
\def\d {\partial\,}
\def\ep{\varepsilon} 
\def\Int {\mbox{\rm int\,}} 

\def\diam {\mbox{\rm diam}\,}

\begin{document}

\title[Weakly convex sets and the modulus of nonconvexity]{Weakly convex sets and modulus of
nonconvexity}

\author{Maxim V. Balashov and Du\v san Repov\v{s}}


\address{Department of Higher Mathematics, Moscow Institute of Physics and Technology, Institutski str. 9,
Dolgoprudny, Moscow region, Russia 141700. balashov@mail.mipt.ru}
\address{Faculty of Mathematics and Physics, and Faculty of Education, University of Ljubljana, Jadranska 19, Ljubljana, Slovenia 1000.
dusan.repovs@guest.arnes.si}

\keywords{Weak convexity, modulus of convexity, modulus of
nonconvexity, proximal smoothness, splitting problem, set-valued
mapping, uniformly continuous selection, uniform convexity.}

\subjclass[2010]{Primary: 54C60, 54C65,52A07.  Secondary: 46A55, 52A01.}

\begin{abstract}
We consider a definition of a weakly convex set which is a
generalization of the notion of a weakly convex set in the sense
of Vial and a proximally smooth set in the sense of Clarke, from
the case of the Hilbert space to a class of Banach spaces with the
modulus of convexity of the second order. Using the new
definition of the weakly convex set with the given modulus of
nonconvexity we prove a new retraction theorem and we obtain new
results about continuity of the intersection of two continuous
set-valued mappings (one of which has nonconvex images) and new
affirmative solutions of the splitting problem for selections. We
also investigate relationship between the new definition and the
definition of a proximally smooth set and a smooth set.

\end{abstract}

\date{\today}
\maketitle

\section{Introduction}

\def\i {\mbox{\rm int}\,}

We begin by some definitions for a Banach space $(E, \|\cdot\|)$
over $\R$. Let $B_{r}(a)=\{ x\in E\ |\ \| x-a\| \le r\}$. Let
$\cl A$ denote the closure and $\Int A$ the interior of the subset
$A\subset E$. The \it diameter \rm of the subset $A\subset E$ is
defined as $\diam A = \sup\limits_{x,y\in A} \| x-y\|$. The \it
distance \rm from the point $x\in E$ to the subset $A\subset E$
is defined as $\varrho (x,A)=\inf\limits_{a\in A}\| x-a\|$. For
a
subset $A\subset E$, let $U_{d}(A)$ be the \it open
$d$-neighborhood \rm of $A$, i.e.
$$
U_{d}(A)=\{ x\in E\ |\ \varrho (x,A)<d\}.
$$
The \it Hausdorff distance \rm between two subsets $A,B\subset\E$
is defined as follows
$$
h(A,B) = \max \left\{ \sup_{a\in A}\ \varrho (a,B) , \quad
\sup_{b\in B}\ \varrho (b,A) \right\}.
$$
We denote the convex hull of the set $A$ by $\co A$.

\begin{Def} (\cite{Balashov+Ivanov09}, \cite{Ivanov})\label{StrConvSegment}
Let $x_{0},x_{1}\in\E$, $\|x_{1}-x_{0}\|\le 2d$. The set
$$
D_{d}(x_{0},x_{1}) = \Bigcap_{a\in\E:\ \{x_{0},x_{1}\}\subset
B_{d}(a)} B_{d}(a)
$$
is called a \it strongly convex segment of radius $d$\rm, and the
set
$$
D^{o}_{d}(x_{0},x_{1}) = D_{d}(x_{0},x_{1})\setminus
\{x_{0},x_{1}\}.
$$
is called a \it strongly convex segment of radius $d$ without
extreme points\rm.
\end{Def}

\begin{Def}\label{VialWeakConv} (Vial \cite{Vial}, see Figure 1).
A  subset $A$ of a normed space is called \it weakly convex (in
the sense of Vial) \rm with constant $R>0$, if for any pair of
points $x_{0},x_{1}\in A$ such that $0<\|x_{1}-x_{0}\|<2R$ the
set $ A\bigcap D^{o}_{R}(x_{0},x_{1})$ is nonempty.
\end{Def}

\begin{Def}\label{Prox} (Clarke et al \cite{Clarke1}, \cite{Clarke}).
A subset $A$ of a normed space $E$ is called \it proximally
smooth \rm with constant $d>0$, if the distance function
$x\to\varrho (x,A)$ is Frechet differentiable on the tube
$U_{d}(A)\backslash A$.
\end{Def}

\begin{Def}\label{modulus} (Polyak \cite{Polyak})
Let $E$ be a Banach space and let a subset $A\subset E$ be convex
and closed. \it The modulus of convexity \rm $\delta_{A}:\
[0,\diam A)\to [0,+\infty) $ is the function defined by
$$
\delta_{A}(\ep) = \sup\left\{ \delta\ge 0\ \left|\
B_{\delta}\left( \frac{x_{1}+x_{2}}{2}\right)\right.\subset A,\
\forall x_{1},x_{2}\in A:\ \| x_{1}-x_{2}\|=\ep \right\}.
$$\rm
\end{Def}

It is obvious that $\delta_{A}(0)=0$.

\begin{Def}\label{RM} (Polyak \cite{Polyak}, see Figure 2).
Let $E$ be a Banach space and let a subset $A\subset E$ be convex
and closed. If the modulus of convexity $\delta_{A}(\ep )$ is
strictly positive for all $\ep\in (0,\diam A)$, then we call the
set $A$ \it uniformly convex (with modulus
$\delta_{A}(\cdot)$).\rm
\end{Def}

We proved in \cite{Balashov+Repovs2} that every uniformly convex
set $A\ne E$ is bounded and if the Banach space  $E$ contains a
nonsingleton uniformly convex set $A\ne E$ then it admits a
uniformly convex equivalent norm. We also proved that the
function $\ep\to \delta_{A}(\ep)/\ep$ is increasing (see also
\cite[Lemma 1.e.8]{Lindestrauss+tzafriri}), and for any uniformly
convex set $A\ne E$ there exists a constant $C>0$ such that
$\delta_{A}(\ep)\le C\ep^{2}$ \cite{Balashov+Repovs2}.

Let $\delta_{E}$ be the modulus of convexity for the Banach space
$E$, i.e. the modulus of convexity for the closed unit ball in
$E$.

\begin{Def}\label{modulus-w}
Let $E$ be a Banach space. Let a subset $A\subset E$ be closed and
$d\in (0,\diam A)$. The \it modulus of nonconvexity \rm
$\gamma_{A}:\ [0,d)\to [0,+\infty) $ is defined as
$$
\gamma_{A}(\ep) = \inf\left\{ \gamma > 0\ \left|\
B_{\gamma}\left( \frac{x_{1}+x_{2}}{2}\right)\right.\cap
A\ne\emptyset,\ \forall x_{1},x_{2}\in A:\ \| x_{1}-x_{2}\|\le\ep
\right\}
$$
and $\gamma_{A}(0)=0$. \rm
\end{Def}

It is easy to see that the modulus of nonconvexity is a
nondecreasing function. Besides, we shall further suppose that the
modulus of nonconvexity is continuous from the right. Otherwise
we shall redefine the modulus by continuity from the right.

\begin{Def}\label{RM-w} (see Figure 3).
Let $E$ be a Banach space, and let a subset $A\subset E$ be
closed. We shall call the set $A$ \it weakly convex with modulus
of nonconvexity $\gamma_{A}(\ep)$, \rm  $\ep\in [0,d)$ ($d\le
\diam A$), if the modulus of nonconvexity $\gamma_{A}$ satisfies
the inequality
$$
 0\le\gamma_{A}(\ep
)<\frac{\ep}{2},\qquad \forall \ep\in [0,d).
$$
\rm
\end{Def}


\begin{figure}[!htb]
\begin{center}
\includegraphics[scale=0.6]{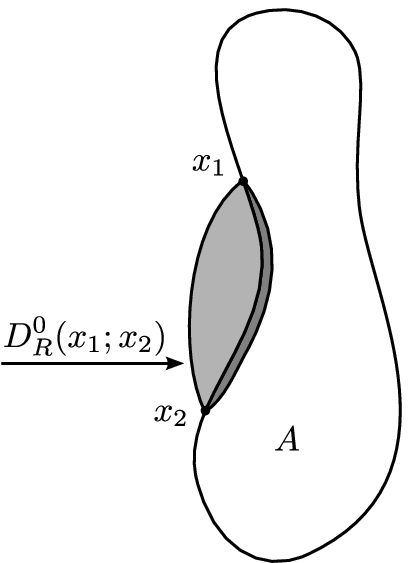}\qquad\qquad\qquad\qquad\includegraphics[scale=0.6]{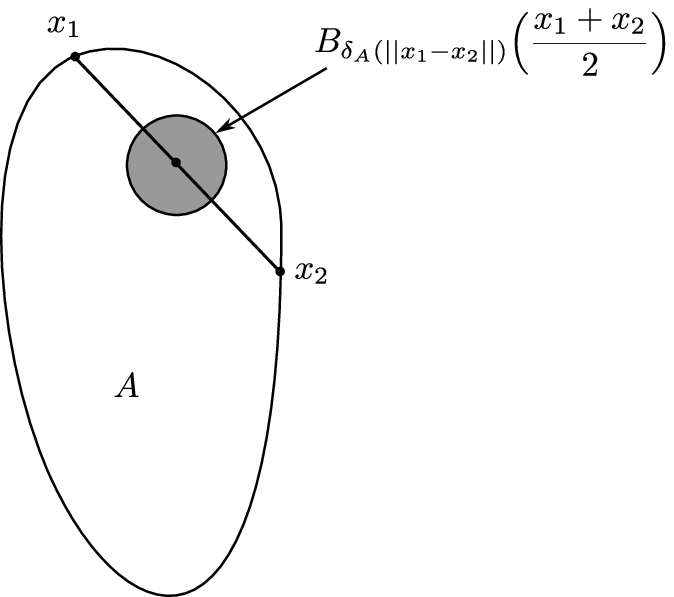}
\qquad \qquad\includegraphics[scale=0.5]{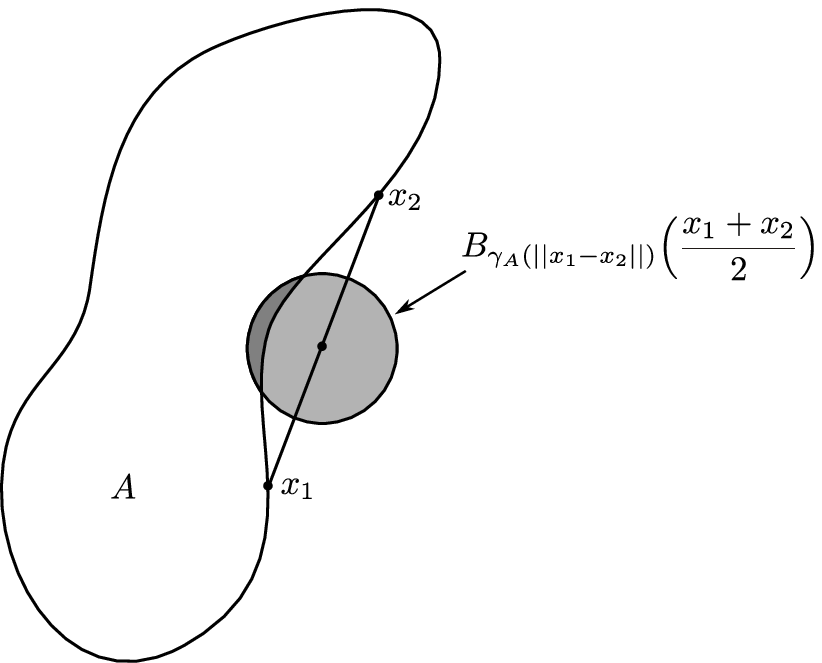}
\end{center}
\vspace{2mm}
Figure 1\qquad\qquad\qquad\qquad\qquad Figure
2\qquad\qquad\qquad\qquad\qquad Figure 3\qquad\qquad \hfill
\end{figure}

It is obvious that the equality $\gamma_{A}(\ep) = 0$ for all
$\ep\in [0,\diam A)$ means (for the closed set $A$) convexity of
the set $A$.

Hereafter the text "weakly convex" means weakly convex in the
sense of Definition 1.\ref{RM-w}.

\Example Let $E=\H$ be the Hilbert space and $\delta_{\H}(\ep) =
1-\sqrt{1-\frac{\ep^{2}}{4}}$ be the modulus of convexity of $\H$.
A weakly convex subset $A\subset \H$  with modulus
$\gamma_{A}(\ep)=d\delta_{\H}(\ep/d)$, $\ep\in [0,d)$, $d>0$, is
weakly convex in the sense of Vial with constant $d$ and
proximally smooth with constant  $d$ (see
\cite{Zlateva,Clarke1,Clarke,Ivanov,Vial}, in particular
\cite{Balashov+Ivanov09}). These three properties are equivalent
in the Hilbert space.

The relationship between weak convexity in the sense of Vial and
proximal smoothness of a set in a Banach space is much more
complicated (see \cite{Balashov+Ivanov09} for details).

The next lemma is a simple consequence of similarity.

\begin{Lm} (\cite[Lemma 2.7.1]{Polovinkin+Balashov})\label{0}
Let a space $E$ be uniformly convex with modulus $\delta_{E}$.
Then for all $x,y\in B_{1}(0)$, such that $\| x-y \|=\ep>0$, and for any
$\beta\in (0,\frac12 ]$ the following inequality holds
$$
B_{2\beta\delta_{E}(\ep)}((1-\beta)x+\beta y)\subset B_{1}(0).
$$
\end{Lm}

\begin{Lm}\label{1}
Let a space $E$ be uniformly convex with modulus $\delta_{E}$.
Then for any $\varepsilon,\eta$ such that
$0<\ep/2<\eta<\ep<2$ the
following inequality holds
$$
\frac{\delta_{E}(\eta)}{\eta}\le
\frac{\delta_{E}(\ep)}{\ep}-2\frac{\ep-\eta}{\ep\cdot\eta}\delta_{E}(r(\ep)),
$$
where $r(\ep)=\frac14 \left(
\frac{\ep}{2}-\delta_{E}(\ep)\right)$.
\end{Lm}

By the Day-Nordlander theorem \cite{Diestel}, $\delta_{E}(\ep)\le
\frac{\ep^{2}}{4}<\ep/2$ for all $\ep\in (0,2)$. Hence $r(\ep)>0$
for all $\ep\in (0,2)$.

\proof Let's fix $\ep\in (0,2)$, $\alpha\in
(0,\frac{\ep}{4}-\frac12\delta_{E}(\ep))$ and $\lambda\in (0,1)$.
Choose points $x_{1},x_{2}\in \d B_{1}(0)$, such that $\|
x_{1}-x_{2}\|=\ep$ and $\delta_{E} (\ep)+\alpha>\delta$, where
$\delta = \sup\{ {r\ge 0}  \ |\ B_{r}(z)\subset B_{1}(0)\}$ and
$z=\frac12 (x_{1}+x_{2})$.

For any natural number $k$ we define the point $a_{k}\in \d
B_{1}(0)$ with the property $\| a_{k}-z\|\le \delta+\frac1k$. Let
$y^{k}_{i}$ be the homothetic image of the point $x_{i}$ under
the homothety with center $a_{k}$ and coefficient $\lambda$,
$i=1,2$; let $z_{k}$ be the homothetic image of the point $z$
under the homothety with center $a_{k}$ and coefficient $\lambda$.

By construction, $\| y^{k}_{1}-y^{k}_{2}\|=\lambda \ep$ and $\|
z_{k}-a_{k}\|\le \lambda \delta+\lambda\frac1k$. By the triangle
inequality and by the property of chosen points $x_{i}$, $a_{k}$
we have $\| x_{i} - a_{k}\|\ge \| \frac{x_{1}-x_{2}}{2}\|-\|
a_{k}-z\|\ge
\frac14\left(\frac{\ep}{2}-\delta_{A}(\ep)\right)=r(\ep)>0$ for
$i=1,2$ and sufficiently large $k$.
 Let $\beta= \min\{ \ll, 1-\ll\}\in (0,\frac12]$. By Lemma 1.\ref{0}
 we have
$B_{2\beta\delta_{E}(r(\ep))}(y^{k}_{i})\subset B_{1}(0)$,
$i=1,2$. Hence
$$
\delta_{E}(\ll\ep)\le \left\|\frac{y^{k}_{1}+y^{k}_{2}}{2}
-a_{k}\right\|-2\beta\delta_{E}(r(\ep))= \|
z_{k}-a_{k}\|-2\beta\delta_{E}(r(\ep)) \le \lambda\delta_{E}(\ep
)+\lambda\alpha+\lambda\frac1k -2\beta\delta_{E}(r(\ep)).
$$
Letting $\alpha\to+0$, $k\to\infty$, we obtain
$$
\delta_{E}(\ll\ep)\le \lambda\delta_{E}(\ep
)-2\beta\delta_{E}(r(\ep)).
$$
The desired estimate appears if we put $\ll=\eta/\ep$. \qed

One of the important
motivations for consideration of
weakly convex sets in
the sense of Definition 1.\ref{RM-w} is given by the next theorem.

\begin{Th}\label{unique-proj}
Let a space $E$ be uniformly convex with modulus $\delta_{E}$,
$d>0$. Let $A\subset E$ be a weakly convex set with modulus of
nonconvexity $\gamma_{A}$, and suppose that function
$d\delta_{E}(\ep/d)-\gamma_{A}(\ep)$ is
positive for all $\ep\in
(0,\min\{ 2d,\diam A\})$. Then for any point $x\in U_{d}(A)$ the
set
$$
P_{A}x = \{ a\in A\ |\ \| x-a\| = \varrho (x,A)\}
$$
is a singleton.
\end{Th}

\proof (1).{\it Nonemptiness of $P_{A}x$}.
Let's fix $x\in
U_{d}(A)\backslash A$. Let points  $a_{k}\in A$ be such that $\|
x-a_{k}\|\to \varrho (x,A)$. Define nonnegative numbers
$\ep_{k}=\| x-a_{k}\| - \varrho (x,A)$.

Suppose that the sequence $\{ a_{k}\}_{k=1}^{\infty}$ has no
converging subsequence. Without loss of generality we may assume
that there exists a number $\ep_{0}>0$ such that for any natural $
k,m$ the following inequality holds: $\|
a_{k}-a_{m}\|\ge\ep_{0}$. By the definition of $\ep_{k}$,
$\ep_{m}$ we have
$$
\max\{ \| x-a_{k}\|, \| x-a_{m}\|\}\le \ep_{k}+\ep_{m}+\varrho
(x,A).
$$
Let $\varrho = \varrho (x,A)$. Then
$$ \left\|
x-\frac{a_{k}+a_{m}}{2}\right\|\le  \varrho
+\ep_{k}+\ep_{m}-(\varrho+\ep_{k}+\ep_{m} )\delta_{E}(\|
a_{k}-a_{m}\|/(\varrho+\ep_{k}+\ep_{m})).
$$
Due to the weak convexity of the set $A$ for any $\alpha>0$ there
exists
$$
a_{km}\in B_{\gamma_{A}(\|
a_{k}-a_{m}\|)+\alpha}\left(\frac{a_{k}+a_{m}}{2}\right)\cap A.
$$
Hence $\| x-a_{km}\|\le  \left\|
x-\frac{a_{k}+a_{m}}{2}\right\|+\gamma_{A}(\|a_{k}-a_{m} \|
)+\alpha\le$
$$
\le\varrho +\ep_{k}+\ep_{m}+\alpha-(\varrho
+\ep_{k}+\ep_{m})\delta_{E}(\| a_{k}-a_{m}\|/(\varrho
+\ep_{k}+\ep_{m})) +\gamma_{A}(\|a_{k}-a_{m} \| ).
$$
Let's choose $d_{1}\in (\frac12 d,d)$ and a sequence
$\alpha_{k}>0$, $\alpha_{k}\to 0$, such that for all sufficiently
large $k,m$ the inequality $\varrho
+\ep_{k}+\ep_{m}+\alpha_{k}<d_{1}<d$ holds. Then by Lemma
1.\ref{1}
$$
(\varrho +\ep_{k}+\ep_{m})\delta_{E}(\| a_{k}-a_{m}\|/(\varrho
+\ep_{k}+\ep_{m}))\ge d_{1}\delta_{E}\left( \frac{\|
a_{k}-a_{m}\|}{d_{1}}\right)
$$
and we  have the estimate $\| x-a_{km}\|\le
\varrho+\ep_{k}+\ep_{m}+\alpha_{k}-$
$$
-d\delta_{E}\left( \frac{\|
a_{k}-a_{m}\|}{d}\right)+\gamma_{A}(\| a_{k}-a_{m}\|) - \left(
d_{1}\delta_{E}\left( \frac{\| a_{k}-a_{m}\|}{d_{1}}\right) -
d\delta_{E}\left( \frac{\| a_{k}-a_{m}\|}{d}\right) \right)\le
$$
$$
\varrho+\ep_{k}+\ep_{m}+\alpha_{k} - \left( d_{1}\delta_{E}\left( \frac{\|
a_{k}-a_{m}\|}{d_{1}}\right) - d\delta_{E}\left( \frac{\|
a_{k}-a_{m}\|}{d}\right) \right).
$$
By Lemma 1.\ref{1} it follows that
$$
d_{1}\delta_{E}\left( \frac{\| a_{k}-a_{m}\|}{d_{1}}\right) -
d\delta_{E}\left( \frac{\| a_{k}-a_{m}\|}{d}\right) \ge
2(d-d_{1})\delta_{E}(r(\| a_{k}-a_{m}\|/d_{1})).
$$
From the inequalities $\ep_{0}\le \| a_{k}-a_{m}\|<2d_{1}$ and
 $r(\ep)=\frac14 \left(
\frac{\ep}{2}-\delta_{E}(\ep)\right)\ge \frac14 \left(
\frac{\ep}{2}-\frac{\ep^{2}}{4} \right)>0$, it follows that for
all $k,m$ the value $\delta_{E}(r(\| a_{k}-a_{m}\|/d_{1}))$ is
bounded from below by a positive constant $c>0$. Hence for
sufficiently large $k,m$ (when
$\ep_{k}+\ep_{m}+\alpha_{k}<2(d-d_{1})c$), $\| x-a_{km}\|<
\varrho$. Contradiction. Therefore, the sequence $\{
a_{k}\}_{k=1}^{\infty}$ has a converging subsequence and
$P_{A}x\ne\emptyset$.

(2).{\it The set $P_{A}x$ is a singleton}.
The proof is similar to the
step 1.
 If $\varrho (x,A)=\|
x-a_{i}\|$, $i=1,2$, $a_{1},a_{2}\in A$, then we have
$$
\left\| x-\frac{a_{1}+a_{2}}{2}\right\|\le \varrho (x,A)-\varrho
(x,A)\delta_{E}(\| a_{1}-a_{2}\|/\varrho (x,A)),
$$
and for all $\alpha>0$
$$
\exists a\in B_{\gamma_{A}(\|
a_{1}-a_{2}\|)+\alpha}\left(\frac{a_{1}+a_{2}}{2}\right)\cap A.
$$
Now by choosing $0<\alpha<\varrho (x,A)\delta_{E}(\|
a_{1}-a_{2}\|/\varrho (x,A))-\gamma_{A}(\| a_{1}-a_{2}\|)$, we
obtain that $\| x-a\|\le \left\| x-\frac{a_{1}+a_{2}}{2}
\right\|+\gamma_{A}(\| a_{1}-a_{2}\|)+\alpha\le$
$$
 \le\varrho (x,A)-
\varrho (x,A)\delta_{E}(\| a_{1}-a_{2}\|/\varrho (x,A))+\gamma_{A}(\| a_{1}-a_{2}\|)+\alpha
<\varrho (x,A).
$$
 Contradiction.\qed

By Theorem 1.\ref{unique-proj} and the results from
\cite{Balashov+Ivanov09} it follows that if the space $E$ is
additionally uniformly smooth then each weakly convex set with
the modulus $\gamma_{A}$ (for which
$d\delta_{E}(\ep/d)-\gamma_{A}(\ep)>0$) is proximally smooth with
constant $d>0$. We note that $d$ is not the largest possible
constant for the proximal smoothness of the set $A$.

It's easy to see that the proximal smoothness with constant $d>0$
implies the weak convexity. Suppose that (for simplicity) the
subset $A\subset E$ is compact in the strong topology of the
Banach space  $E$ and proximally smooth with constant $d>0$. Then
the set $A$ is weakly convex with some modulus of nonconvexity
$\gamma_{A}(\ep)$, $\ep\in (0, \min\{2d,\diam A\})$. Indeed, the
compactness of the set $A$ implies that the values of modulus from
the Definition 1.\ref{modulus-w} are achieved for every $\ep\in
(0, \min\{2d,\diam A\})$
  at some points $a_{1},a_{2}\in A$, $\| a_{1}-a_{2}\|\le\ep$. This means that for the point $x=\frac12 (a_{1}+a_{2})$
  we have  $\varrho (x,A)=\gamma_{A}(\ep)\ge 0$. Using inequality $\varrho (x,A)\le \frac12\ep$,
  we obtain from the estimate $\frac12\ep<d$ and from proximal smoothness
  of the set  $A$ (see \cite[Theorem 2.4]{Balashov+Ivanov09}) that
  the set $P_{A}x$ is a singleton and $\varrho
(x,A)=\gamma_{A}(\ep)<\frac12\ep$.

\section{The order of function $\gamma_{A}$}

Before further considerations we shall make some remarks.
Consider for simplicity a set $A$ on the Euclidean plane. Let
the boundary $\d A$ be a smooth closed curve  $x=x(s)$, $y=y(s)$,
where $s$ is the natural parameter. Suppose that the curve $\d A$
contains no straight segments. In this case the radius of 
curvature of $\d A$ at the point  $(x(s),y(s))$ equals
$R(s)=(x'^{2}(s)+y'^{2}(s))^{3/2}/|x''(s)y'(s)-y''(s)x'(s)|$. If
the radius $R(s)$ is finite and positive at the point
$(x(s),y(s))$
 (and this takes place for a.e. values of parameter $s$),
then the curve at the neighborhood of the point  $(x(s),y(s))$ is
similar to the circle of radius $R(s)$. 

If additionally, the set
$A$ is not locally convex at the point $(x(s),y(s))$ (i.e. for
any $r>0$ the set $A\cap B_{r}((x(s),y(s)))$ is nonconvex), then
for a small $\ep>0$ the function $\gamma_{A}(\ep)$ has the order
no smaller than $\ep^{2}$ (more precisely, $\gamma_{A}(\ep)\ge
R^{2}(s)-\sqrt{R^{2}(s)-\frac{\ep^{2}}{4}}$).

We shall show that the situation above is typical: if the set $A$
is nonconvex, then the modulus of nonconvexity for $A$ satisfies
the estimate $\gamma_{A}(\ep)\ge \mbox{\rm Const}\cdot \ep^{2}$.
As we have mentioned above the modulus of nonconvexity for convex
set $A$ equals zero.

If the subset $A$ of the Banach space $E$ is a symmetric \it
cavern \rm (i.e. $A=\cl(E\backslash B)$, where $B$ is a closed
convex bounded symmetric body), then $\gamma_{A} (\ep)\ge
\mbox{\rm Const}\cdot \ep^{2}$. The proof follows by the fact that
in this case the function $\gamma_{A}(\ep)$ has the same order as
the the function
$$ \sigma_{E,B}(\ep)=\sup\left\{ 1-\frac{\| x+y\|_{B} }{2}\ \left|\
x,y\in \d B,\ \| x-y\|_{B}\le\ep\right.\right\},
$$
introduced in \cite{Banas1}. Here $\|\cdot\|_{B}$ is the norm in
the space $E$ with the unit ball $B$. In \cite{Banas1} and
\cite{Banas2} the inequality $1-\sqrt{1-\frac{\ep^{2}}{4}}\le
\sigma_{E,B}(\ep)$ was proved. In fact, it is the "dual" of the
Day-Nordlander theorem. It can be proved similarly as the
Day-Nordlander theorem (see \cite[\S 3, pp. 60--62]{Diestel}
for details). The proof is the same except that instead of
function $\delta_{X}(\ep) = \inf_{\varphi}\Delta (\ep,\varphi)$
one should consider on page 62 the function $\sigma_{X}(\ep) =
\sup_{\varphi}\Delta (\ep,\varphi)$.

Let the subset $A\subset E$ from a Banach space $E$ be a \it
cavern, \rm i.e. $A=\cl (E\backslash B)$ where the set $B\subset
E$ is a closed convex and bounded body, $0\in \i B$. We shall
estimate the value of $\gamma_{A}(\ep)$.

For any closed convex and bounded set $B\subset E$, $0\in\i B$,
we define the \it Minkowski \rm function
$$
\mu_{B}(x) = \inf\{ t>0\ |\ x\in tB\},\quad \forall x\in E.
$$
For any  bounded set $C\subset E$ we define \it $B$-diameter \rm
of the set $C$ as follows
$$
\diam_{B}C=\sup\limits_{x,y\in C}\mu_{B}(x-y).
$$
For any closed convex bounded sets $A,B,C\subset E$ we  define
the modulus
$$
\sigma_{C}^{A,B}(\ep) = \inf\left\{ \sigma\ge 0\ |\ \left(\sigma
A+\frac{x+y}{2}\right) \cap (E\backslash \i C)\ne\emptyset,\quad
\forall x,y\in A:\ \mu_{B}(x-y)\le\ep\right\}
$$
and the modulus
\begin{equation}\label{BanasMod}
\sigma_{C}(\ep) = \sigma_{C}^{B_{1}(0),B_{1}(0)}(\ep).
\end{equation}

Moduli $\sigma_{C}^{A,B}$ and $\sigma_{C}$ generalize the
definition from \cite{Banas1} to arbitrary convex sets.

It is obvious from the definition of $\sigma_{C}$ that if $C$ is
a convex body then we have for all admissible $\ep>0$ for the set
$A=\cl (E\backslash C)=\cl (E\backslash \i C)$
\begin{equation}\label{g>s}
\gamma_{A}(\ep) = \sigma_{C}(\ep).
\end{equation}

The next lemmas are direct consequences of the definition of
$\sigma_{C}^{A,B}$.

\begin{Lm}\label{s-1}
For any bounded closed convex bodies $A,B,C\subset E$ and $t>0$,
the following holds:
$$
 \sigma_{C}^{tA,B}(\ep) =
\frac{1}{t}\sigma_{C}^{A,B}(\ep),\quad \forall\ep\in
(0,\diam_{B}C);\leqno(1)
$$
$$
\sigma_{C}^{A,tB}(\ep) = \sigma_{C}^{A,B}(t\ep),\quad
\forall\ep\in \left(0,\frac1t\,\diam_{B}C\right);\quad\mbox{\rm
and}\leqno(2)
$$
$$
\sigma_{tC}^{A,B}(\ep) =
t\sigma_{C}^{A,B}\left(\frac{\ep}{t}\right),\quad \forall\ep\in
\left(0,t\,\diam_{B}C\right).\leqno(3)
$$
\end{Lm}

\begin{Lm}\label{s-2}
For any bounded closed convex bodies $A',B',A,B,C\subset E$ and
$\ep\in
(0,\diam_{B}C)$,\\
(1) if $A'\subset A$ then $\sigma_{C}^{A',B}(\ep)\ge
\sigma_{C}^{A,B}(\ep)$; and\\
(2) if $B'\subset B$ then $\sigma_{C}^{A,B'}(\ep)\le
\sigma_{C}^{A,B}(\ep)$.
\end{Lm}

\begin{Th}\label{g-2}
Suppose that the subset $A\subset E$ of a Banach space $E$ is a
cavern. Let $B_{r}(0)\subset \cl (E\backslash A)\subset
B_{R}(0)$. Then for all $\ep\in (0,2r)$ we have
$$
\gamma_{A}(\ep)\ge \frac{\ep^{2}}{8R^{2}}r.
$$
\end{Th}

\proof Let $B=\cl (E\backslash A)$ be a closed convex body. Using
Lemmas 2.\ref{s-1} and 2.\ref{s-2} we get
$$
\sigma_{B}(\ep)=\sigma_{B}^{B_{1}(0),B_{1}(0)}(\ep) =
r\sigma_{B}^{B_{r}(0),B_{R}(0)}\left(\frac{\ep}{R}\right)\ge
r\sigma_{B}^{B,B}\left(\frac{\ep}{R}\right).
$$

Using the result of Bana\'{s} \cite{Banas1} we obtain 
$$
\sigma_{B}^{B,B}\left(\frac{\ep}{R}\right)\ge
\sigma_{\H}\left(\frac{\ep}{R}\right) =
 1-\sqrt{1-\frac{\ep^{2}}{4R^{2}}}\ge
\frac{\ep^{2}}{8R^{2}}.
$$

By invoking formula (2.\ref{g>s}) we complete the proof.\qed

Of course, a weakly convex set is not necessarily a cavern. But if
such set $A$ is connected and nonconvex, then it has "cavern-like"
parts, and hence $\gamma_{A}(\ep)\ge\mbox{\rm Const}\cdot
\ep^{2}$.

Hereafter all Banach spaces will have the modulus of convexity of
the second order at zero and will contain weakly convex sets with
modulus of nonconvexity of the second order, too. There are many
such spaces besides Hilbert spaces, for example $l_{p}$, $p\in
(1,2)$ (see \cite{Balashov+Ivanov09}, \cite{Zlateva},
\cite{Diestel} for details).

We shall define a special condition.

\begin{Def}\label{M}
Let $\delta$ be the modulus of convexity for some closed convex
set $A$ and $\gamma$ the modulus of nonconvexity for some closed
weakly convex set $B$. We shall say that \it condition (i) \rm is
valid for the moduli $\delta$ and $\gamma$ if \\
(1) for all $s\in [0,s_{0}]$ there exists a solution $t=t_{s}>s$
of the equation $\delta (t-s)-\gamma (t)=0$,\\
(2) the function $\delta (t-s)-\gamma (t)$ is positive and
increasing for $t>t_{s}$ and\\
(3) there exists a solution $t(s)$ of the equation $\delta
(t-s)-\gamma (t)=\frac{s}{2}$.

In the case when $\delta (t) = d\delta_{E}(t/d)$ is the modulus
of convexity for the ball $B_{d}(0)$ in the Banach space $E$ we
shall define the solution $t(s)$ as $t_{E}(s)$.
\end{Def}

\begin{Remark}\label{zam-1}
Definition 2.\ref{M} has a technical character (it is useful for
further proofs) and it is not so exotic. Suppose that for
sufficiently small $t>0$ our moduli have the second order at zero
and are defined by formulae $\delta(t)=c_{1}t^{2}+o(t^{2})$,
$t\to+0$, $\gamma(t)=c_{2}t^{2}+o(t^{2})$, $t\to+0$, and
$c_{1}>c_{2}>0$. Then for sufficiently small numbers $s>0$ the
function $t\to \delta(t-s)- \gamma(t)$ is positive and increasing,
and $t(s)\asymp\sqrt{s}$, $s\to +0$.
\end{Remark}

\begin{Remark}\label{zam-11}
Suppose that a Banach space $E$ has modulus of convexity
$\delta_{E}$ of the second order and  a closed subset $A\subset E$
is weakly convex with the modulus of nonconvexity $\gamma_{A}$ of
the second order, too. Taking into account that for all
$0<d<d_{1}$ $d\delta_{E}(\ep/d)\ge d_{1}\delta_{E}(\ep/d_{1})$,
$\forall \ep\in (0,2d)$ (see \cite[Lemma 2.1]{Balashov+Repovs2}),
and $d\delta_{E}(\ep/d)\asymp\frac{\ep^{2}}{d}$, $\ep\to+0$, we
can conclude that there exists a number $d>0$ such that
$d\delta_{E}(\ep/d)>\gamma_{A}(\ep)$, $\forall \ep\in (0,2d)$. If
additionally the space $E$ is smooth then by Theorem
1.\ref{unique-proj} and \cite[Theorem 2.4]{Balashov+Ivanov09} we
obtain that the set $A$ is proximally smooth with constant $d$.
\end{Remark}

\section{Properties of weakly convex sets}

\begin{Th}\label{svyaznost2} Let $d>0$. Let a subset $A$ of a Banach space $E$  be weakly convex with modulus
$\gamma_{A}(\ep)$, $\ep\in [0,d)$ and the  subset $B\subset E$ be
uniformly convex with modulus $\delta_{B}(\ep)$, $\ep\in [0,\diam
B)$ and $\diam B<d$. Let $\delta_{B}(\ep) > \gamma_{A}(\ep)$ for
all $\ep\in [0, \diam B)$. Then the set $A\cap B$, if nonempty,
is weakly convex with modulus $\gamma_{A\cap B}(\ep
)\le\gamma_{A}(\ep)$, $\ep\in [0,\diam A\cap B)$ and connected.
\end{Th}

\proof The weak convexity of the intersection and the estimate
for the modulus follows by definitions.

Suppose that the set $A\cap B$ is not connected.  This means that
there exist two nonempty closed disjoint sets $A_{1}\subset A\cap
B$ and $A_{2}= ( A\cap B)\backslash A_{1}$. Choose $k=1$ and
points $a_{1}\in A_{1}$, $b_{1}\in A_{2}$.

Due to weak convexity of the set  $A\cap B$ there exists a point
$w\in \frac12 (a_{k}+b_{k})+(\gamma_{A}(\|
a_{k}-b_{k}\|)+\alpha_{k})B_{1}(0)$, $w\in A$. The numbers
$\alpha_{k}$ are chosen by the conditions  $\alpha_{k}\to 0 $
 and $0<\alpha_{k}<\frac12 (\frac12\| a_{k}-b_{k}\|-\gamma_{A}(\| a_{k}-b_{k}\|))$.

One of the inclusions  $w\in A_{1}$ or $w\in A_{2}$ is true. If
$w\in A_{1}$, then denote $a_{k+1}=w$, $b_{k+1}=b_{k}$. If $w\in
A_{2}$, then denote $a_{k+1}=a_{k}$, $b_{k+1}=w$. In this way we
build the sequences  $\{ a_{k}\}_{k=1}^{\infty}\subset A_{1}$, $\{
b_{k}\}_{k=1}^{\infty}\subset A_{2}$.

Let $l_{k}=\| a_{k}-b_{k}\|$. Then $0\le l_{k+1}\le \frac12
l_{k}+\gamma_{A}(l_{k})+\alpha_{k}<l_{k}$. Hence $l_{k}\to l\ge
0$. Taking the limit $k\to \infty$ and using continuity of the
function $\gamma_{A}$ from the right we get $\frac12 l\le
\gamma_{A}(l)$. It follows from Definitions 1.\ref{modulus-w} and
1.\ref{RM-w} that $l=0$. Therefore, $\| a_{k}-b_{k}\|\to 0$.

We proved in \cite{Balashov+Repovs2} that for any uniformly
convex set $B$ there exists a number $c>0$ such that the modulus
of convexity for the set  $B$ can be estimated as follows
$\delta_{B}(\ep)\le c\ep^{2}$, $\ep\in (0,\diam B)$.

It follows from the construction of points $a_{k+1}$ that
$a_{k+1}=a_{k}$, or
$$
\| a_{k+1}-a_{k}\|\le \frac12 \|a_{k}-b_{k}\|+\gamma_{A}(\| a_{k}-b_{k}\|)+\alpha_{k}\le   \frac34 l_{k}+\frac{c}{2}l_{k}^{2}.
$$
In the latter case $b_{k+1}=b_{k}$ and
$$
l_{k+1} =\| a_{k+1}-b_{k}\|\le \frac34 l_{k}+\frac{c}{2}l_{k}^{2} = \left( \frac34+\frac{c}{2}l_{k}\right)l_{k}\le \frac45 l_{k},
$$
for all $k>k_{0}$. Thus there exists a number  $d>0$, such that
$l_{k}\le d\left(\frac45\right)^{k}$. It follows from the estimate
$$
\| a_{k+1}-a_{k}\|\le \frac{3d}{4}\left( \frac{4}{5}\right)^{k}+\frac{c}{2}d \left( \frac{4}{5}\right)^{2k}\le K\left( \frac{4}{5}\right)^{k},
$$
which is valid for sufficiently large $k$, that for such $k$ and
$m>k$
$$
\| a_{m}-a_{k}\|=\sum\limits_{n=k}^{m-1}\| a_{n+1}-a_{n}\|\le \sum\limits_{n=k}^{m-1} K\left( \frac{4}{5}\right)^{n} \le
5K\left( \frac{4}{5}\right)^{k},
$$
the latter means that the sequence $\{a_{k}\}$ is fundamental.
Since $\| a_{k}-b_{k}\|\to 0$ thus the sequence $\{ b_{k}\}$ is
also fundamental. By the closedness of the sets $A_{1}$ and
$A_{2}$ and from the condition $\| a_{k}-b_{k}\|\to 0$ we
conclude that $a_{k}\to x\in A_{1}$, $b_{k}\to x\in A_{2}$. Hence
$A_{1}\cap A_{2}\ne \emptyset$.\qed

For any closed subset $A\subset E$, a point $x\in U_{d}(A)$ and a
number $s>0$ we define the \it set-valued projection \rm
$$
P_{A}(x,s)=\{ a\in A\ |\ \| x-a\|\le \varrho (x,A)+s\}.
$$
It follows by definition that $P_{A}(x,s)\ne\emptyset$ for all
$s>0$. Apart from this, under conditions of Theorem
1.\ref{unique-proj}, $P_{A}(x,0)=P_{A}x$ is a singleton.

\begin{Th}\label{pr-meat}
Let a Banach space $E$ be uniformly convex with modulus
$\delta_{E}$. Let a subset $A\subset E$ be weakly convex with
modulus $\gamma_{A}(\ep)$, $\ep\in [0,2d)$. Let
$d\delta_{E}(\ep/d)> \gamma_{A}(\ep)$ for all $\ep\in (0,2d)$.
Suppose that the condition (i) from Definition 2.\ref{M} is
satisfied. Then
$$
P_{A}(x,s)\subset B_{t_{E}(s)}(P_{A}x).
$$
\end{Th}

\proof Let $a=P_{A}x$ and $b\in P_{A}(x,s)$, $s\le s_{0}$. Let's
define the point $y\in [b,x]$ by the condition  $\| x-y\|=\varrho
(x,A)=\| x-a\|$. Let $w=\frac{a+b}{2}$, $z=\frac{a+y}{2}$.

It follows from the inequality $\| b-y\|\le s$ that $\| w-z\|\le
s/2$. In the triangle $bya$ we see that $\| y-a\|\ge \| a-b\|-s$.

Let $\varrho =\varrho (x,A)$. Note that $\| a-b\|-s\le \| y-a \|<
2\varrho< 2d$. If the inequality $\varrho\delta_{E}((\|
a-b\|-s)/\varrho)-\gamma_{A}(\| a-b\| )>\frac{s}{2}$ holds, then
for some $\alpha>0$ we have $\varrho\delta_{E}((\|
a-b\|-s)/\varrho)-\gamma_{A}(\| a-b\| )>\frac{s}{2}+\alpha$.
Using the inequality $\varrho\delta_{E}((\| a-b\|
-s)/\varrho)<\varrho\delta_{E}((\| a-y\|)/\varrho)$, we obtain
$$
a_{0}\in B_{\gamma_{A} (\| a-b\|)+\alpha}(w)\cap A,\qquad  B_{\gamma_{A}
(\| a-b\|)+\alpha}(w)\subset \Int B_{\varrho\delta_{E}((\| a-b\|
-s)/\varrho)}(z).
$$

Hence $\| a_{0}-x\|\le \| a_{0}-w\|+\| w-z\|+\| z-x\|\le$
$$
\le
\gamma_{A} (\| a-b\| )+\alpha+\frac{s}{2}+\varrho (x,A)-\varrho
(x,A)\delta_{E}((\| a-b\| -s)/\varrho (x,A))<\varrho (x,A).
$$
This contradiction shows that $d\delta_{E} ((\|
a-b\|-s)/d)-\gamma_{A}(\| a-b\| )< \varrho\delta_{E}((\|
a-b\|-s)/\varrho)-\gamma_{A}(\| a-b\| )\le\frac{s}{2}$ and by the
conditions of the theorem, $\| a-b\|\le t_{E}(s)$. The point
$b\in P_{A}(x,s)$ was arbitrary and the theorem is thus
proved.\qed

\begin{Corol}\label{sl-1}
Under the assumptions of Theorem 3.\ref{pr-meat}, the projection
$P_{A}x$ uniformly continuously depends on $x$. More precisely, if
$\| x_{1}-x_{2}\| <s_{0}$ and $x_{1},x_{2}\in U_{d}(A)$, then $\|
P_{A}x_{1}-P_{A}x_{2}\|\le t_{E}(\| x_{1}-x_{2}\|)$. Moreover,
$t_{E}(s)\asymp\sqrt{s}$, $s\to+0$.
\end{Corol}

\begin{Th}\label{svyaznost1}  Suppose that
the assumptions of Theorem 3.\ref{pr-meat} hold and $d_{1}\in
(0,d)$.  Then for any point $x\in E$ the set $A\cap
B_{d_{1}}(x)$, if nonempty, is weakly convex with modulus
$\gamma_{A\cap B_{d_{1}}(x)}(\ep )\le\gamma_{A}(\ep)$, $\ep\in
[0,\diam A\cap B_{d_{1}}(x))$, and path connected.
\end{Th}

\proof Weak convexity of the intersection follows from the
definitions.

Fix any pair of points $x,y\in A$ such that $0<\|x-y\|<2d_{1}$.
For any number $t\in[0;1]$ we denote $z_{t}=(1-t)x+t y$. The map
 $z\mapsto P_{A}z$ is single-valued and continuous (Corollary 3.\ref{sl-1}) on the set
$U_{d_{1}}(A)$, hence it is single-valued and continuous on the
set
 $U'=U_{d_{1}}(A)\bigcup A$. Since $z_{t}\in U'$ for all $t\in[0;1]$
 there is
a unique point $a(t)$ with $\{a(t)\}= P_{A}z_{t}$. The function
$a:[0;1]\to A$ is continuous and defines the desired curve
$\Gamma=\{a(t):\ t\in[0;1]\}$ which connects points $x$ and
$y$.\qed


\begin{figure}[!htb]
\centering
\includegraphics[scale=1.0]{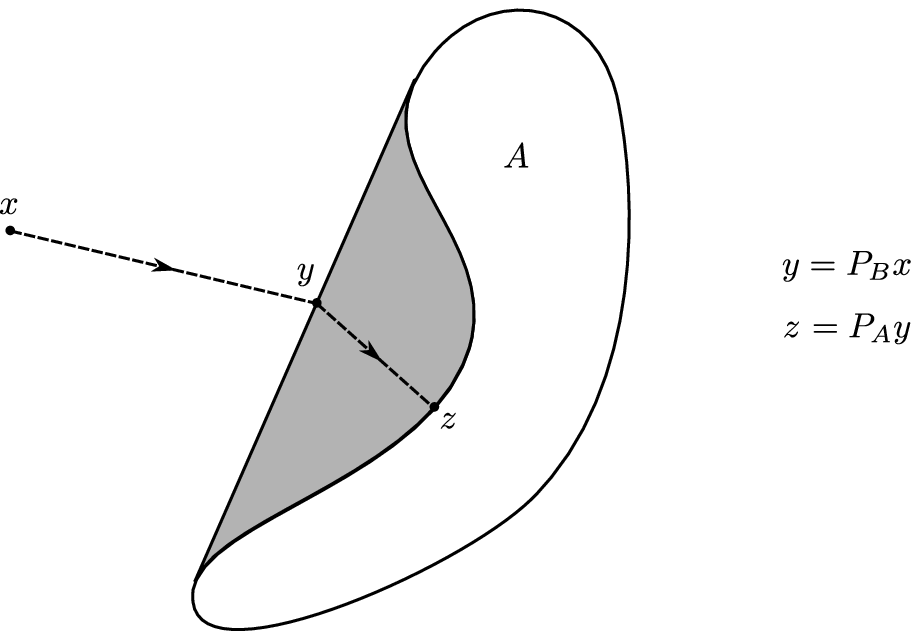}
\vspace{2mm}
\\Figure 4: Scheme of retraction
\end{figure}

\begin{Th}\label{retract}
Let $E$ be a uniformly convex space with modulus  $\delta_{E}$.
Let $A\subset E$ be a weakly convex set with modulus of
nonconvexity $\gamma_{A}(\ep)$, $\ep\in [0,\diam A)$. Let
$A\subset B_{r}(a)$, $2r<d$ and $\gamma_{A}(\ep)<
d\delta_{E}(\ep/d)$ for all $\ep\in [0,\diam A)$. Suppose that
the assumptions of Theorem 3.\ref{pr-meat} hold. Then the set $A$
is a continuous retract of $E$.
\end{Th}

\proof  Let $x\in E\backslash A$.
 Let $B=\cl\co A\subset B_{r}(a)$,
 $y=P_{B}x$. Due to the uniform convexity of the space  $E$ the
metric projection on the set  $B$ is continuous. We observe that
this projection is uniformly continuous (see \cite{Lindenstrauss}
and \cite[Example 3.2]{Balashov+Repovs2}) on the balls.

Since $y=y(x)\in B\subset B_{r}(a)$ we have $\varrho (y,A)\le
2r<d$. By Theorem 1.\ref{unique-proj} and Corollary 3.\ref{sl-1}
there exists a unique metric projection $z=P_{A}y$ which
uniformly continuously depends on  $y$. Therefore,
$z(x)=P_{A}(P_{B}x)$
 is the desired retraction, see Figure 4.\qed

\begin{Remark}\label{sl-2}
We remark that function $z(x)$ from Theorem 3.\ref{retract} is
uniformly continuous on the balls.

Let us also mention that Theorem 3.\ref{retract} remains valid in
any uniformly convex and smooth Banach space for any proximally
smooth set $A$ with constant $d$ and $A\subset B_{r}(a)$, $d<2r$.
Instead of Theorem 1.\ref{unique-proj} and Corollary 3.\ref{sl-1}
one must use the results from \cite[Theorem
2.4]{Balashov+Ivanov09}.
\end{Remark}

\begin{Th}\label{fin-din-cont}
Let  $(T,\rho)$ be a metric space. Let $F_{1},F_{2}:(T,\rho)\to
2^{\R^{n}}$ be set-valued mappings, continuous in the Hausdorff
metric. Suppose that for a point $t_{0}\in T $ the set
$F_{1}(t_{0})$ is uniformly convex with modulus  $\delta (\ep)$,
and the set $F_{2}(t_{0})$ is weakly convex with modulus  $\gamma
(\ep)$. Let $\gamma (\ep)<\delta (\ep)$ for all $\ep<\min\{\diam
F_{1}(t_{0}),\ \diam F_{2}(t_{0})\}$. Let $H(t)=F_{1}(t)\cap
F_{2}(t)\ne\emptyset$ for all $t\in T$. Then the mapping $H(t)$ is
continuous at the point $t=t_{0}$ in the Hausdorff metric.
\end{Th}

\proof It follows from the uniform convexity of the set
$F_{1}(t_{0})$ that it is bounded (\cite{PL}, \cite[Theorem
2.1]{Balashov+Repovs2}). Due to the continuity in the Hausdorff
metric we conclude that there exists a number $\delta>0$ such that
the set $\cl \bigcup\limits_{\rho (t,t_{0})<\delta}F(t)$ is
compact. By the Closed Graph Theorem \cite{Aubin} the set-valued
mapping $H(t)$ is upper semicontinuous at the point $t=t_{0}$,
i.e.
$$ \forall\ep>0\ \exists \delta>0 \ \forall t:\ \rho
(t,t_{0})<\delta\qquad H(t)\subset H(t_{0})+B_{\ep}(0),
$$
or
\begin{equation}\label{semi-up}
\lim\sup\limits_{t\to t_{0}}H(t)\subset H(t_{0}).
\end{equation}
If the set $H(t_{0})$ is a singleton then the continuity of the
set-valued mapping $H$ at the point $t_{0}$ follows by its upper
semicontinuity. Next we shall assume that the set $H(t_{0})$
consists of more than one point.

Suppose that lower semicontinuity fails, i.e.
 that
 $$
H(t_{0})\not\subset \lim\inf\limits_{t\to t_{0}}H(t).
$$
Thus there exist a number $\ep_{0}>0$
 and points $t_{k}\in T$, $t_{k}\to t_{0}$, such that
 $$
H(t_{0})\not\subset H(t_{k})+B_{\ep_{0}}(0),\qquad\mbox{\rm for
any natural}\,\, k.
 $$
For any $k$ there exists a point $h_{k}\in H(t_{0})$ with
$$
h_{k}\not\in H(t_{k})+B_{\ep_{0}}(0).
$$
Since the set $H(t_{0})$ is compact, thus without loss of
generality we may assume that $h_{k}\to h_{0}\in H(t_{0})$ and
\begin{equation}\label{notsemi-low}
h_{0}\not\in H(t_{k})+B_{\ep_{0}/2}(0),\qquad \mbox{\rm for any
natural}\,\, k.
\end{equation}

Let us define the set
$H_{0}=\lim\sup\limits_{k\to\infty}H(t_{k})\subset H(t_{0})$. By
construction $h_{0}\in H(t_{0})\backslash H_{0}$, hence $H_{0}\ne
H(t_{0})$.

Let $x_{0}\in P_{H_{0}}h_{0}$. For $h_{0}\in F_{1}(t_{0})\cap
F_{2}(t_{0})$ and $x_{0}\in H_{0}$ put $l=\| h_{0}-x_{0}\|>0$.

Using uniform convexity of $F_{1}(t_{0})$ we get
$$
B_{\delta (l)}\left( \frac{x_{0}+h_{0}}{2}\right)\subset
F_{1}(t_{0}).
$$
Due to the weak convexity of the set $F_{2}(t_{0})$ and finite
dimension of images of the mapping $F_{2}$ we can find a point
$$
f\in B_{\gamma (l)}\left( \frac{x_{0}+h_{0}}{2}\right)\cap
F_{2}(t_{0}).
$$
By continuity of the map $F_{2}$ there exist points $f_{k}\in
F_{2}(t_{k})$ with $f_{k}\to f$. Besides, for $\ep = (\delta
(l)-\gamma (l))/3$ we can find a natural number $k_{0}$, such
that for all $k>k_{0}$ the following holds:
\begin{equation}\label{h1}
f_{k}\in B_{\delta (l)-\ep}\left(
\frac{x_{0}+h_{0}}{2}\right)\subset F_{1}(t_{k}),
\end{equation}
and
\begin{equation}\label{h2}
f_{k}\in B_{\gamma (l)+\ep}\left(
\frac{x_{0}+h_{0}}{2}\right)\cap F_{2}(t_{k}).
\end{equation}
By the formulae (3.\ref{h1}) and (3.\ref{h2}) it follows that
$f_{k}\in H(t_{k})$ for all $k>k_{0}$. Hence
$f=\lim\limits_{k\to\infty}f_{k}\in H_{0}$.

At the same time
$$
\| h_{0}-f\|\le \left\| h_{0}-\frac{x_{0}+h_{0}}{2}\right\|+
\left\| f-\frac{x_{0}+h_{0}}{2}\right\|\le \frac12\|
x_{0}-h_{0}\|+\gamma (\| x_{0}-h_{0}\|)<\| x_{0}-h_{0}\|.
$$
This contradicts with the inclusion $x_{0}\in P_{H_{0}}h_{0}$.
Thus, $H(t)$ is lower semicontinuous at the point $t=t_{0}$.\qed

Let $F:(T,\rho)\to 2^{(E,\|\cdot\|)}$ be a set-valued mapping. If
for any $t\in T$ the set $F(t)$ is uniformly convex with modulus
$\delta_{F(t)}(\ep)\ge\delta(\ep)$, $\ep\in [0,\diam F(t))$, and
$\delta $ is an increasing function, then we shall say that the
set-valued mapping $F$ is uniformly convex with modulus $\delta$.

If for any $t\in T$ the set $F(t)$ is weakly convex with modulus
of nonconvexity $\gamma_{F(t)}(\ep)\le\gamma(\ep)$, $\ep\in
[0,\diam F(t))$, $\gamma (0)=0$, $\gamma (\ep)<\frac{\ep}{2}$ for
admissible $\ep
> 0$ and function $\gamma$ is continuous from the right and nondecreasing then we shall say that the set-valued mapping
$F$ is uniformly weakly convex with modulus $\gamma$.

\begin{Def}\label{M1}
Let a set-valued mapping  $F_{1}:(T,\rho)\to 2^{(E,\|\cdot\|)}$ be
uniformly convex with modulus $\delta$ and a set-valued mapping
$F_{2}:(T,\rho)\to 2^{(E,\|\cdot\|)}$ uniformly weakly convex
with modulus $\gamma$. We shall say that \it condition (ii) \rm
is valid for the moduli
$\delta$ and $\gamma$ if \\
(1) for all $s\in [0,s_{0}]$ there exists a solution $t=t_{s}>s$
of the equation $\delta (t-s)-\gamma (t)=0$,\\
(2) the function $\delta (t-s)-\gamma (t)$ is positive and
increasing for $t>t_{s}$,\\
(3) there exists a solution $t(s)$ of the equation $\delta
(t-s)-\gamma (t)=\frac{s}{2}$.

\end{Def}

It follows from the results of the second paragraph that
condition (ii) is possible only if moduli $\delta$ and $\gamma$
are of the second order at zero.

We say that set-valued mapping $F$  is uniformly continuous with
modulus of continuity $\omega\ge 0$ if for any $t_{1},t_{2}\in T$
the inequality $h(F(t_{1}),F(t_{2}))\le \omega (\rho
(t_{1},t_{2}))$ holds.

\begin{Th}\label{intersection}  Let $F_{1},F_{2}:(T,\rho)\to 2^{(E,\|\cdot\|)}$.
Let the values of $F_{2}$ be uniformly convex with modulus
$\delta (\ep)$. Let the values of $F_{1}$ be uniformly weakly
convex with modulus $\gamma (\ep)$. Suppose that set-valued
mapping $F_{i}$ is uniformly continuous with modulus
$\omega_{i}$, $i=1,2$. Let the condition (ii) holds.

Let $H(t)=F_{1}(t)\cap F_{2}(t)\ne\emptyset$ for all $t\in T$ and
suppose that for some $M>0$ the inclusion $\bigcup\limits_{t\in
T}H(t)\subset B_{M}(0)$ holds. Then
\begin{equation}\label{cap} h (H(t_{1}), H(t_{2}))\le
 \left\{\begin{array}{l}
2\omega_{1}+3\omega_{2}+t\left(\frac{\omega_{1}+\omega_{2}}{2}\right),
\quad \frac{\omega_{1}+\omega_{2}}{2}<s_{0},\\
\frac{\omega_{1}+\omega_{2}}{s_{0}}M,\quad
\frac{\omega_{1}+\omega_{2}}{2}\ge s_{0}.
\end{array}\right.
\end{equation}
\end{Th}
Here $\omega_{i}=\omega_{i}(\rho (t_{1},t_{2}))$, $i=1,2$.

\proof Let $c_{1}\in H(t_{1})$. We shall show that for any number
$\ll$, which is strictly larger than the right side of the formula
(3.\ref{cap}), there exists a point $c_{2}\in H(t_{2})$ with $\|
c_{1}-c_{2}\|\le \ll$. This will prove the theorem.

Fix $d\in H(t_{2})$. If $\omega_{1}+\omega_{2}\ge 2s_{0}$, then,
taking $c_{2}=d$, we obtain that
$$
h(H(t_{1}),H(t_{2}))\le \| c_{1}-c_{2}\|\le 2M\le \frac{\omega_{1}+\omega_{2}}{s_{0}}M.
$$

Suppose that $\omega_{1}+\omega_{2}< 2s_{0}$. Fix $k>1$, such
that inequality $k\omega_{1}+k^{2}\omega_{2}< 2s_{0}$ holds. For
the point $c_{1}\in H(t_{1}) = F_{1}(t_{1})\cap F_{2}(t_{1})$ we
can find the point $b\in F_{2}(t_{2})$ such that $\| b-c_{1}
\|\le k\omega_{2}$.

Fix the point $b_{\pi}\in F_{1}(t_{2})$ which satisfies the
condition $\| b-b_{\pi}\|\le k\cdot\varrho (b,F_{1}(t_{2}))$.
Invoking the inequality $\varrho (b,F_{1}(t_{1}))\le \|
b-c_{1}\|\le k\omega_{2}$ we get the following estimate
$$
\| b-b_{\pi}\|\le k\varrho (b,F_{1}(t_{2}))\le
kh(F_{1}(t_{1}),F_{1}(t_{2}))+k\varrho (b,F_{1}(t_{1}))\le
k\omega_{1}+k\| b-c_{1}\|\le k\omega_{1}+k^{2}\omega_{2}.
$$

Define the point $a\in [d,b]\cap H(t_{2})$ as the one which is
nearest to the point $b$. The set $[d,b]\cap H(t_{2})$ is nonempty
because it contains the point $d$. Put $n=1$, $a_{1}=a$.

Consider the following cases:

{\bf (1)}  $\delta (\| a_{n}-b\|)> \gamma (\|
a_{n}-b_{\pi}\|)+\frac12 \| b-b_{\pi}\|$ or

{\bf (2)}  $\delta (\| a_{n}-b\|)\le \gamma (\|
a_{n}-b_{\pi}\|)+\frac12 \| b-b_{\pi}\|$.

If the case (1) takes place then we choose $\alpha_{n}=$
$$
\min\left\{ \frac{1}{n}, \frac12 \left( \delta (\| a_{n}-b\|)- \gamma (\|
a_{n}-b_{\pi}\|)-\frac12 \| b-b_{\pi}\|\right) , \frac12 \left( \frac{\| a_{n}-b_{\pi}\|}{2}- \gamma (\|
a_{n}-b_{\pi}\|)\right)\right\}>0.
$$
By the uniform weak convexity of $F_{1}$ with the modulus
$\gamma$ there exists a point
\begin{equation}\label{f*2}
 w\in B_{\gamma (\| a_{n}-b_{\pi}\|)+\alpha_{n}}\left(
\frac{a_{n}+b_{\pi}}{2}\right)\cap F_{1}(t_{2})\subset B_{\delta
(\| a_{n}-b\|)}\left( \frac{a_{n}+b}{2}\right)\subset
F_{2}(t_{2}),
\end{equation}
and

$$
\| b_{\pi}-w\|\le \left\|
b_{\pi}-\frac{a_{n}+b_{\pi}}{2}\right\|+\left\|
\frac{a_{n}+b_{\pi}}{2}-w\right\|\le \frac12 \|
a_{n}-b_{\pi}\|+\gamma (\| a_{n}-b_{\pi}\|)+\alpha_{n} <\| a_{n}-b_{\pi}\|.
$$
Now we put $n=n+1$, $a_{n}=w$ and again consider cases (1) or (2).

If the case (2) does not take place for all natural  $n$, then we
obtain from the construction of the points $\{a_{n}\}$  that the
sequence $l_{n}=\| a_{n}-b_{\pi}\|$ satisfies the condition $0\le
l_{n+1}\le \frac{l_{n}}{2}+\gamma (l_{n})+\alpha_{n}<l_{n}$. It
follows by the Weierstrass theorem that the sequence $l_{n}$
converges to some number $l\ge 0$ from the right. Using the
continuity of the function $\gamma$ from the right and taking the
limit we deduce that $\frac{l}{2}\le \gamma (l)$. The latter is
possible only in the case $l=0$ (see the definition of $\gamma$).

Thus, if for all $n\in\N$ the case (2) does not take place then
$H(t_{2})\ni a_{n}\to b_{\pi}$, i.e. $b_{\pi}\in H(t_{2})$. Taking
$c_{2}=b_{\pi}$ we have
$$
\| c_{1}-c_{2}\|=\| c_{1}-b_{\pi}\|\le \| c_{1}-b\| +\|b-b_{\pi}
\|\le k\omega_{2}+k\omega_{1}+k^{2}\omega_{2}.
$$
The number $k>1$ was arbitrary, hence
$$
h(H(t_{1}),H(t_{2}))\le \omega_{1}+2\omega_{2}.
$$

Suppose that for some $n\in\N$ the case (2) occurs and $\|
a_{n}-b_{\pi}\| > \| b-b_{\pi}\|$. Taking into account that $\|
a_{n}-b\|>\| a_{n}-b_{\pi}\|-\| b-b_{\pi}\|$, we conclude from
the inequality of the case (2), that
$$
\delta (\| a_{n}-b_{\pi}\|-\| b-b_{\pi}\|)-\gamma (\|
a_{n}-b_{\pi}\|)\le \frac12 \| b-b_{\pi}\|.
$$
From the condition (ii) of the theorem we get
$$
\| a_{n}-b_{\pi}\|\le t\left(\frac12 \| b-b_{\pi}\|\right)\le
t\left(\frac{k\omega_{1}+k^{2}\omega_{2}}{2}\right).
$$

By choosing  $c_{2}=a_{n}$ we obtain
$$
h(H(t_{1}),H(t_{2}))\le \| c_{1}-c_{2}\|\le \| c_{1}-b\|+\|
b-b_{\pi}\|+\| a_{n}-b_{\pi}\| \le
k\omega_{2}+k\omega_{1}+k^{2}\omega_{2}+t\left(\frac{k\omega_{1}+k^{2}\omega_{2}}{2}
\right).
$$
By taking the limit $k\to 1+0$, we finally prove the theorem. The
case $\| a_{n}-b_{\pi}\| \le \| b-b_{\pi}\|$, which follows from
the last estimate and from the inequality $\| b-b_{\pi}\|\le
k\omega_{1}+k^{2}\omega_{2}$, also gives formula
(3.\ref{cap}).\qed

\begin{Remark}\label{g1}
For convex valued mapping $F_{2}$ a similar result was proved in
\cite[Theorem 3.1]{Balashov+Repovs2}.
\end{Remark}

\begin{Remark}\label{g10}
If additionally the conditions of the Theorem 3.\ref{svyaznost2}
hold for sets $F_{1}(t)$ and $F_{2}(t)$, then the values of the
map $H$ in Theorem 3.\ref{intersection} are connected.
\end{Remark}

\begin{Remark}\label{g2}
In our case the moduli $\delta$ and $\gamma$ are of the second
order and we have that $t(s)$ is of the order $\sqrt{s}$ when
$s\to 0$. For the Hilbert space this result was proved by Ivanov
\cite{Ivanov}.
\end{Remark}

\section{Application to selection problems}

\begin{Th}\label{selection} Let $E$ be a uniformly convex Banach space with modulus  $\delta_{E}$.
Let $\Phi\subset E$ be a collection of sets with weakly convex
images with modulus of nonconvexity  $\gamma (\ep)$ (in the sense
of Section 3), and suppose that all sets from $\Phi$ are contained
in some ball. Let $d>0$. Let the condition (ii) holds for moduli
$d\delta_{E} (t/d)$ and $\gamma (t)$. Let $d\delta_{E}
(t/d)>\gamma (t)$ for all admissible $t>0$. Suppose that any set
$H\in \Phi$ is contained in (each in its own) ball of radius
$r>0$ and $2r<d$.

Then the collection $\Phi$ has a uniformly continuous selection,
i.e. there exists a uniformly continuous in the Hausdorff metric
function   $s:\Phi\to E$  such that for all $H\in \Phi$ we have
$s(H)\in H$.
\end{Th}

\proof Without loss of generality we shall assume that the sets
from the family  $\Phi$
 are contained in the ball  $B_{R}(0)$  and for any $H\in \Phi$ we have $\varrho (0,\cl\co H)\ge r_{1}>0$.
 Consider $\Psi=\{ \cl\co H\ |\ H\in \Phi\}$. The metric projection of zero
  $y(H)=P_{\cl \mbox{\footnotesize{\rm co}}\, H}0$ on the sets from $\Psi$ is a uniformly continuous selection
  defined of $\Psi$.
   We have proved  in \cite[Lemma 3.1]{Balashov+Repovs2} that for any $H_{1},H_{2}\in \Phi$
  (taking into account the condition $h(\cl \co H_{1},\cl \co H_{1})\le h(H_{1},H_{1})$)
 $$
\| y(H_{1})-y(H_{2})\|\le 2h (H_{1},H_{2})+f_{E}(h(H_{1},H_{2})),
 $$
where
$$
f_{E}(t)=\left\{\begin{array}{l} \delta^{-1}(t/2), \quad
t<2\Delta_{E},\\
\frac{Rt}{\Delta_{E}},\quad t\ge 2\Delta_{E.}
\end{array}\right.
$$
Here $\delta ( \ep)=R\delta_{E}(\ep/R)$, $\Delta_{E}=\delta
(2r_{1})$.

Let $y=y(H)$. From $\varrho (y,H)\le 2r <d $ using the Theorem
1.\ref{unique-proj} we conclude that there exists a unique metric
projection $z(H)=P_{H}y$.

If $2h (H_{1},H_{2})+f_{E}(h(H_{1},H_{2}))<(d-2r)/2$, then, by
defining $y_{i}=y(H_{i})$, $z_{i}=z(H_{i})$, $i=1,2$, we get $\|
y_{1}-y_{2}\|< (d-2r)/2$.

Consider a metric subspace $T$ of the metric space $((E,\Phi),(\|
\cdot, \cdot\|+h(\cdot,\cdot)))$. Elements of $T$ are pairs
$(x,H)\in (E,\Phi)$ such that $\varrho (x,H)<d$. Consider the
set-valued mappings $F_{1}(x,H)=H$, $F_{2}(x,H)=B_{\varrho
(x,H)}(x)$ from $T$ into $E$. The set-valued mapping $F_{1}$ is
uniformly weakly convex with modulus $\gamma$ and uniformly
continuous with modulus $\omega_{1}(t)=t$. The set-valued mapping
$F_{2}$ is uniformly convex with modulus $d\delta_{E} (\ep/d)$.

For points $(y_{i},H_{i})$, $i=1,2$, we have
$$h(F_{2}(y_{1},H_{1}),F_{2}(y_{2},H_{2}))=\| y_{1}-y_{2}\|+|\varrho (y_{1},H_{1})-\varrho (y_{2},H_{2})|\le$$
$$
\le\| y_{1}-y_{2}\|+|\varrho (y_{1},H_{1})-\varrho (y_{1},H_{2})|+|\varrho (y_{1},H_{2})-\varrho (y_{2},H_{2})|,
$$
$|\varrho (y_{1},H_{1})-\varrho (y_{1},H_{2})|\le
h(H_{1},H_{2})$, and from the condition $\| y_{1}-y_{2}\|\le
(d-2r)/2$ we obtain that
$$
\varrho (y_{1},H_{2})\le \| y_{1}-y_{2}\|+\varrho (y_{2},H_{2})\le (d-2r)/2 + 2r = (d+2r)/2<d.
$$
Put $z_{12}\in H_{2}$: $\| y_{1}-z_{12}\|=\varrho (y_{1},H_{2})$.
Using Corollary  3.\ref{sl-1} we get $|\varrho
(y_{1},H_{2})-\varrho (y_{2},H_{2})|\le \| y_{1}-y_{2}\|+\|
z_{2}-z_{12}\|\le \| y_{1}-y_{2}\| + t_{E}(\| y_{1}-y_{2}\|)$.

Thus in the case $2h
(H_{1},H_{2})+f_{E}(h(H_{1},H_{2}))<(d-2r)/2$ projections
$z_{i}=F_{1}(y_{i},H_{i})\cap F_{2}(y_{i},H_{i})$ uniformly
continuously depend on sets $H_{i}$, $i=1,2$, by Theorem
3.\ref{intersection}, i.e.
$$
\| z_{1}-z_{2}\|\le \omega (h(H_{1},H_{2})),
$$
where $\omega (h (H_{1},H_{2}))$ is superposition of the function
$2h (H_{1},H_{2})+f_{E}(h(H_{1},H_{2}))$ and the function from the
right side of formula (3.\ref{cap}).

If $2h (H_{1},H_{2})+f_{E}(h(H_{1},H_{2}))\ge (d-2r)/2$, then by
the strict monotonicity (increasing) of the function $f_{E}$,
there exists a number $C>0$, such that $h(H_{1},H_{2})>C$. In this
case
$$
\| z_{1}-z_{2}\|\le 2R\le \frac{2R}{C}h(H_{1},H_{2})
$$
Therefore, $s(H)=z(H)$ is a uniformly continuous selection.\qed

\Example One can apply these results to certain questions about
continuous selections of set-valued mappings \cite{LM-LS},
\cite{Repovs+Semenov}. Let a space $E$ be uniformly convex, a sets
$A,B\subset E$ be such that $B$ is uniformly convex with modulus
$\delta (\ep)$, $\ep\in[0,\diam B)$, and $A$ is weakly convex
with modulus $\gamma (\ep)$, $\ep\in [0,\diam A)$. Let for some
$d>0$ the inequalities $2\diam B<d$ and $\gamma
(\ep)<d\delta_{E}(\ep/d)$ for all $\ep\in [0,\min \{ 2d,\diam A
\})$ hold. Suppose that the condition (i) is valid for pairs
$\delta (\ep)$, $\gamma (\ep)$  and $d\delta_{E} (\ep/d)$,
$\gamma (\ep)$. Then there exist uniformly continuous functions
$a:A+B\to A$ and $b:A+B\to B$ such that for any  $c\in A+B$ we
have $a(c)+b(c)=c$.

\proof By Theorem 3.\ref{intersection} the set-valued mapping
$A+B\ni c\to H(c)=B\cap (c-A)$ is uniformly continuous. By
definition the set $H(c)$ is weakly convex with modulus of
nonconvexity $\gamma $ for all $c\in A+B$. By the boundedness of
the set $B$, all sets $H(c)$, $c\in A+B$, are contained in some
ball. Furthermore for any point $c\in A+B$ each set $H(c)$ is
contained in the ball of radius no larger than  $2\diam B<d$. By
Theorem 4.\ref{selection} there exists a uniformly continuous
selection $b(c)=s(H(c))\in B$, where the function $s(\cdot)$ is
from Theorem 4.\ref{selection}; $a(c)=c-b(c)\in A$.\qed

\section{A class of weakly convex sets}

We shall show that simple smooth closed surfaces of codimension 1
are weakly convex sets. In this section the space $E$ will be an
arbitrary reflexive Banach space.

 We introduce the \it normal cone \rm $N(A,x)$ to
the set $A$ at the point $x\in A$ as follows
$$
N(A,x) = \{ p\in E^{*}\ |\ (p,x-a)\ge -\alpha_{x}(\|
x-a\|)\cdot\| x-a\|\cdot \| p\|,\quad\forall a\in A\},
$$
where the function $\alpha_{x}:[0,\diam A)\to [0,+\infty)$ and
$\lim_{t\to +0}\alpha_{x}(t)=0$.

Let $A\subset E$ be any closed set with the property $\cl\i A =A$
and $x\in A$. Suppose that the set $\d A$ has the following
properties: $\d A$ is path connected, $\forall x\in \d A$
$$
N(\d A,x)\cap \d B^{*}_{1}(0)= \left(N(A,x)\cap \d
B^{*}_{1}(0)\right)\bigcup \left( N(\cl (E\backslash A),x)\cap \d
B^{*}_{1}(0)\right),
$$
where
$$
N(A,x)\cap \d B^{*}_{1}(0) = \{ p\},\quad  N(\cl (E\backslash
A),x)\cap \d B^{*}_{1}(0) = \{-p\},
$$
 and there exists infinitely small at zero function $\alpha
:[0,\diam A)\to [0,+\infty)$ with the property
$$
N(A,x) = \{ p\in E^{*}\ |\ (p,x-a)\ge -\alpha(\| x-a\|)\cdot\|
x-a\|\cdot\| p\|,\quad\forall a\in A\},\quad\forall x\in \d A,
$$
$$
N(\cl (E\backslash A),x) = \{ p\in E^{*}\ |\ (p,x-a)\ge
-\alpha(\| x-a\|)\cdot\| x-a\|\cdot\| p\|,\quad\forall a\in \cl
(E\backslash A)\}, \forall x\in \d A.
$$
  Then we say that the set $\d A$ is a \it smooth
closed surface of codimension 1 with a function of smoothness
$\alpha$. \rm Roughly speaking, smooth closed surface of
codimension 1 is the smooth path connected
 boundary between some  set $A$ and its complementary set $\cl (E\backslash
 A)$.

Let $r>0$. Define for any point $x\in \d A$ and for unit vector
$p\in N(\d A, x)$ the vector $y\in E$ with $\| y\|=(p,y)=1$. We
say that a smooth closed surface $\d A$ is \it simple, \rm if for
any 2-dimensional affine plane $L$, such that $\{x,x+y\}\subset
L$, the intersection $L\cap \d A\cap B_{r}(x)$ is a path
connected planar curve. We call $r>0$ the \it parameter of
simplicity.\rm


\begin{Th}\label{class} Let $E$ be a reflexive Banach space. Suppose that $A\subset E$ is a closed set, $\cl
\i A= A$, and $\d
 A$ is a simple smooth closed surface of codimension 1 with the function of smoothness
 $\alpha$ and the parameter of simplicity $r>0$.
 Then the set $\d A$ is weakly
 convex with the modulus of nonconvexity $\gamma_{A}(\ep)\le
 \ep\left(\alpha (\ep)+\alpha  (\ep/2)\right)$ for all $\ep \in [0,\min\{r, \ep_{0}\})$; where
 $\ep_{0} = \sup\{ t>0\ |\ \alpha (\tau)+\alpha  (\tau/2)<\frac12,\quad \forall
 \tau\in (0,t)\}$.
\end{Th}

\begin{figure}[!htb]
\centering
\includegraphics[scale=1.0]{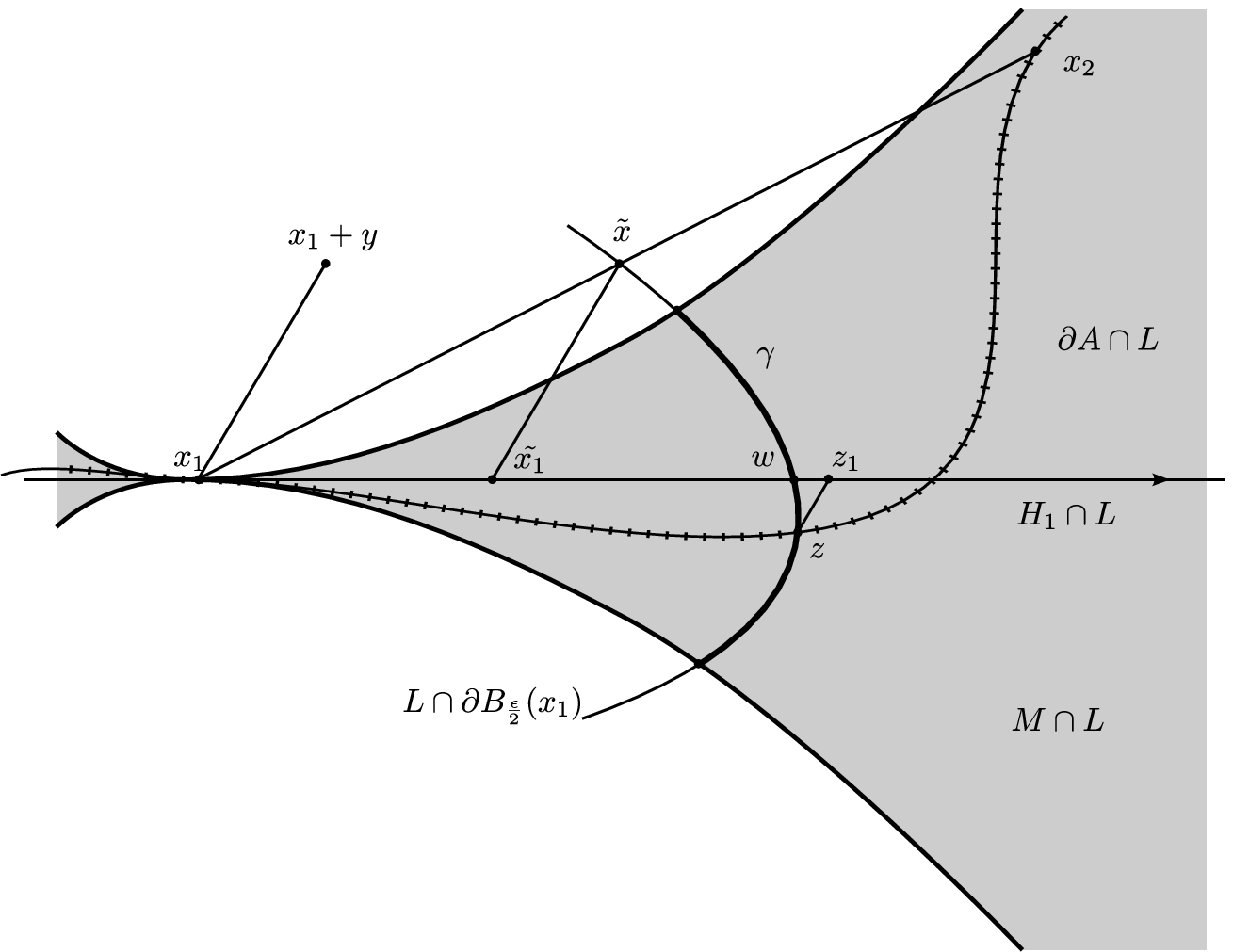}
\vspace{3mm}
\\Figure 5: Proof of Theorem 5.1
\end{figure}

\proof Let $\ep\in (0,\min\{r, \ep_{0}\})$, $x_{1},x_{2}\in \d A$
such that $\| x_{1}-x_{2}\|=\ep$. Let $p_{1}\in N(A,x_{1})\cap \d
B^{*}_{1}(0)$ or $p_{1}\in N(\cl (E\backslash A) ,x_{1})\cap \d
B^{*}_{1}(0)$. Let $H_{1}=\{ x\in E\ |\
(p_{1},x_{1}-x)=0\}=x_{1}+\ker p_{1}$.

Consider Figure 5. Define $M=\{ x\in E\ |\ |(p_{1},x_{1}-x)|\le \alpha
(\| x_{1}-x\|)\cdot\| x_{1}-x\|\}$. Using the reflexivity of the
space $E$ let $y\in E$, $\| y\|=1$, $(p_{1},y)=1$. Let $L=\aff\{
x_{1},x_{1}+y,x_{2}\}$. Let $\tilde x = \frac12 (x_{1}+x_{2}) $,
from the definition of $L$ and from the definition of smooth
closed surface $\varrho (x_{2},H_{1}\cap L) = \varrho
(x_{2},H_{1})\le \alpha (\ep)\ep$, hence $\varrho (\tilde x,
H_{1}\cap L)\le \frac{\ep}{2} \alpha (\ep)$. Let $\tilde x_{1}\in
H_{1}\cap L$ be a point such that $\| \tilde x - \tilde
x_{1}\|=\varrho (\tilde x, H_{1}\cap L)$.

Let $\gamma$ be the connected part of the planar curve $L\cap \d
B_{\ep/2}(x_{1})\cap M$,  which intersects the line $H_{1}\cap L$,
and lies in the same hyperplane with the point $\tilde x$ with
respect to the line $\aff\{ x_{1},x_{1}+y\}$. By the simplicity
of the surface $\d A$ and by the inequality $\ep<r$ there exists
$z\in \gamma\cap\d A$. From the inclusion $z\in M$ we have
$\varrho (z,H_{1}\cap L)\le \alpha (\ep/2)\frac{\ep}{2}$. Let
$z_{1}\in H_{1}\cap L$ be a point such that $\| z -
z_{1}\|=\varrho (z, H_{1}\cap L)$.

Choose the right direction of the line $H_{1}\cap L$ from the
point $x_{1}$ to the point $w=H_{1}\cap L\cap\gamma$.

By the triangle inequality we have from the triangle $x_{1}\tilde
x \tilde x_{1}$
$$
\frac{\ep}{2}-\frac{\ep}{2}\alpha (\ep) \le \| x_{1}-\tilde
x_{1}\|\le \frac{\ep}{2}+\frac{\ep}{2}\alpha (\ep)
$$
and from the triangle $x_{1}zz_{1}$
$$
\frac{\ep}{2}-\frac{\ep}{2}\alpha\left( \frac{\ep}{2}\right) \le
\| x_{1}-z_{1}\|\le \frac{\ep}{2}+\frac{\ep}{2}\alpha\left(
\frac{\ep}{2}\right).
$$
If the point $\tilde x_{1}$ lies to the right of the point $w$
then
$$
\| \tilde x_{1} - w\|=\| x_{1}-\tilde x_{1}\|-\| x_{1}-w\|\le
\frac{\ep}{2}\alpha (\ep).
$$
If the point $\tilde x_{1}$ lies to the left of the point $w$ then
$$
\| \tilde x_{1} - w\|=\| x_{1}-w\|-\| x_{1}-\tilde x_{1}\|\le
\frac{\ep}{2}\alpha (\ep).
$$
In both cases $\| \tilde x_{1}-w\|\le \frac{\ep}{2}\alpha ( \ep)$.

In the same way we obtain that $\| z_{1}-w\|\le
\frac{\ep}{2}\alpha\left( \frac{\ep}{2}\right)$.

From this we deduce that
$$
\| z_{1}-\tilde x_{1}\|\le \| \tilde x_{1}-w\|+\| z_{1}-w\|\le
\frac{\ep}{2}\alpha ( \ep) + \frac{\ep}{2}\alpha\left(
\frac{\ep}{2}\right).
$$
Finally, for any $x_{1},x_{2}\in \d A$, $\| x_{1}-x_{2}\|=\ep$,
$\tilde x = \frac12 (x_{1}+x_{2})$ there exists  $z\in \d A$ such
that
$$
\| \tilde x - z\|\le \| \tilde x - \tilde x_{1}\|+\| z-z_{1}\|+\|
\tilde x_{1}-z_{1} \|\le  \ep\left( \alpha ( \ep) + \alpha\left(
\frac{\ep}{2}\right) \right) <\frac{\ep}{2}.
$$

\qed

We observe that under the conditions of Theorem 5.\ref{class} both
sets $A$ and $\cl (E\backslash A)$ are weakly convex. The proof
easily follows from Theorem 5.\ref{class}.

Let $E$ be a Banach space and a subset $A\subset E$ be closed. We
shall say that unit vector $n\in E$ is a {\it proximall normal
\rm} to the set $A$ at the point $x\in \d A$ if there exists
$r>0$ such that
$$
A\cap \i B_{r}(x+rn)=\emptyset.
$$

\begin{Th}
Let space $E$ be uniformly convex with modulus of convexity of
the second order and uniformly smooth, and let subsets $A\subset
E$ and $\cl (E\backslash A)$ be weakly convex with modulus $\gamma
(\ep)$ of the second order and $\cl \i A=A$. Then  $N(A,x)\cap \d
B_{1}^{*}(0)=\{ p(x)\} $, $N(\cl(E\backslash A),x)\cap \d
B_{1}^{*}(0)=\{ -p(x)\} $,  at any point $x\in\d A$ and $p(x)$
uniformly continuously depends on the point $x\in \d A$.
\end{Th}

\proof By Remark 2.\ref{zam-11} the sets $A$ and $\cl (E\backslash
A)$ are proximally smooth with some parameter $d>0$. Using the
results of Ivanov \cite[Theorem 2]{Ivanov2} we have that for any
point $x\in \d A$ there exists a proximally normal vector  $n(x)$
to the set $A$ at the point $x\in \d A$ (and proximally normal
vector  $-n(x)$ to the set $\cl (E\backslash A)$ at the point
$x\in \d A$) and $n(x)$ uniformly continuously depends on $x$. By
uniform smoothness of the space $E$ for any vector $n\in E$, $\|
n\|=1$, there exists unit vector $p(n)\in E^{*}$ with
$(p(n),n)=1$ and $p(n)$ uniformly continuously depends on $n$
(\cite{Diestel,Lindestrauss+tzafriri}). Again by the smoothness of
the space $E$ we have that $p(x)=p(n(x))\in N(A,x)$ is uniformly
continuous on $x\in \d A$.\qed

\section{Epilogue}

{\bf 1.} We see from the  results above
that in the spaces with
modulus of convexity of the second order the notion of weakly
convex set is very effective.

{\bf 2.} Some of the results can be proved in a more general setting of
uniformly convex Banach spaces: Theorem 3.\ref{retract}, or
Theorems 3.\ref{svyaznost2}, 3.\ref{pr-meat}, and 3.\ref{svyaznost1} 
(see \cite{Balashov+Ivanov09} for details). However, the proofs in
\cite{Balashov+Ivanov09} are much more complicated.

We wish to point out that the  results above
are interesting and
nontrivial even in the finite-dimensional case.

{\bf 3.} We agree with Bana\'{s} that the modulus $\sigma$ from
\cite{Banas1} and \cite{Banas2} is sometimes much more convenient
in applications than the standard modulus of smoothness
\cite{Diestel}, \cite{Lindestrauss+tzafriri}. In fact, the
modulus of nonconvexity is a modification of modulus $\sigma$ from 
\cite{Banas1} for the nonconvex case.

Also, Theorem 5.\ref{class} shows a deep relationship between
weakly convex sets and smooth sets. In conclusion, we formulate
the following

 \begin{Conj}
 Let the
space $E$ be uniformly convex (and smooth?), the subsets
$A\subset E$ and $\cl (E\backslash A)$ closed and weakly convex
with modulus $\gamma (\ep)$ and $\cl \i A=A$. If
$\lim\limits_{\ep\to+0}\gamma (\ep)/\ep=0$, then the unit normal
vector to the set $A$ at the point $x\in\d A$ uniformly
continuously depends on the point $x\in \d A$. 
\end{Conj}

\bigskip

 \section*{Acknowledgements}
\bigskip

  This research was supported by SRA grants P1-0292-0101,
J1-9643-0101, J1-2057-0101, and BI-RU/08-09/001. The first author
was supported by  RFBR grant 10-01-00139-a, ADAP project
"Development of scientific potential of higher school" 2.1.1/500
and project of REA 1.2.1 NK-13P/4.
We thank the referee for several comments and suggestions.

\end{document}